\documentclass[12pt]{article}
\usepackage{amssymb}
\newtheorem{theorem}{Theorem}
\newtheorem{corollary}{Corollary}

\newtheorem{proposition}{Proposition}

\title{Memory Efficient Arithmetic}
\author{Ernie Croot}
\date{\today}

\begin{document}
\maketitle
\begin{abstract}  

In this paper we give an algorithm for finding the $m$th base-$b$ 
digit of a positive integer $n$
($m=1$ is the least significant digit) defined as the final number
in a sequence of integers gotten by multiplying, adding, and subtracting
previous numbers in the sequence (actually, the algorithm finds arbitrarily
precise approximations to $n/b^m \pmod{1}$, which can be used to get this $m$th 
digit whenever the lower $m-1$ digits do not begin with a long run of the digit
$b-1$).  In many cases, this algorithm will require far 
less memory than it takes to write down the base-$b$ digits of $n$, and will have 
a running time (in bit operations) only slightly worse than linear in the number of these base-$b$ digits.

One easy-to-state consequence of the above 
result is that the $m$th base-$10$ digit of 
$2^t$ can be found using $O(t^{2/3} \log^C t)$ bits of memory and 
$O(t \log^C t)$
bit operations, where $C > 0$ is constant.  Of course, if $m = O(t^{2/3})$
then one can do much better (by just computing $n \pmod{b^m}$), 
so the result is only non-trivial when $t = O(m^{3/2})$.

The algorithm we give is highly parallelizable, although if one uses
$M$ processors, to get an $M$-fold reduction in running time, the memory requirements
will increase by a factor of $M$.
\end{abstract}

\section{Introduction}

Suppose that $\alpha$ is a positive real number.  Then, a standard fact regarding
base-$b$ representations is that given any integer $b \geq 2$,
there exists and integer $J$, and a sequence of integers $r_J,r_{J-1},...$, such that
\begin{equation} \label{b_expansion}
\alpha\ =\ \sum_{j \leq J} {r_j b^j},\ 0 \leq r_j \leq b-1,\ r_J \neq 0,
\end{equation}
and we write
$$
\alpha\ =\ (r_J r_{J-1}...r_0.r_{-1}...)_b
$$
to denote this expansion (when $J > 0$).  If we further disallow $r_n = b-1$ for all 
$n \leq N$, for some $N$, then the sequence of $r_i$'s is unique.  
If we make no such restrictions on $r_n$, then the sequence of $r_i$'s
may not be unique, as in the base-$10$ expansion
$$
1\ =\ 0.9999999...
$$

A question which has attracted recent attention (see \cite{bailey})
is the following:
Given a base $b$, and an integer $m \geq 0$, efficiently determine 
a good approximation to 
$$
\nu\ =\ \{ b^m \pi \},
$$
where for a real number $\theta$, the notation $\{\theta\}$ means
$$
\{ \theta\}\ =\ \theta\ -\ \lfloor \theta \rfloor
\ \ ({\rm Note:}\ \ \{\theta\} \equiv \theta \pmod{1},\ 0 \leq \{\theta \} < 1).
$$  
We will say that such an approximation $\gamma$ is a
{\it level-$p$ approximation} if and only if 
$$
|\gamma - \nu|\ <\ {1 \over b^p}.
$$
Now, if the approximation $\gamma$ is sufficiently good, then one can use it 
to determine the digits of $\pi$:  For example, if we take $m = 3$ and $b=10$, 
then we get
$$
\nu\ \equiv\ b^m \pi\ \equiv\ 10^3 \pi\ \equiv\ 
0.59265358979... \pmod{1}.
$$ 
Suppose we had an approximation $\gamma$ to $\nu$
satisfying
$$
|\gamma\ -\ \nu|\ <\ {1 \over 1000};
$$
so, $\gamma = 0.59...$.  Then, the leading base-$10$ digit of $\gamma$
is the same as the fifth base-$10$ digit of $\pi$ (from the left), 
which is $5$.

More generally, a good approximation to $\nu$ gives us
the $(m+2)$nd digit of $\pi$; however, depending on the value of $m$ 
selected, this approximation may need to be extremely close to $\nu$,
in order to determine this digit.  For example, in the above 
instance with $m=3$, if $\gamma = 0.6$, then
$$
|\gamma\ -\ \nu|\ <\ {1 \over 100},
$$
and we note that the leading digit of $\gamma$ is not the same as the 
leading digit of $\nu$.  
\bigskip

In this paper we will describe a method for determining the $m$th
base-$b$ digit of an integer, where $m=1$ corresponds to the 
least significant digit (note that for the digits-of-$\pi$ problem above,
$m=1$ corresponded to the leading digit; so, the $m$th digit is defined 
differently in this context).  As in the problem concerning digits of $\pi$,
this method produces arbitrarily precise approximations to  
$$
\nu\ =\ \left \{ {n \over b^m} \right \},
$$
where $n$ is this integer, which will be defined by a certain type of 
expansion that we will describe below.  Note that the value of $m$ is 
at least $1$ (if $m \leq 0$, then $\nu$ is trivially $0$).  Now, if 
$$
\nu\ =\ (0.r_{-1} r_{-2} \cdots )_b,
$$
then the $m$th digit of $n$ equals $r_{-1}$; and, if $\lambda$ is a
sufficiently good approximation to $\nu$, then the leading digit of $\lambda$
will also equal $r_{-1}$.  

The type of expansion for $n$ we will use is defined as follows:
In \cite{shub}, Smale and Shub say that a {\it computation of length }$L$ 
of a positive integer $n$ is a sequence of integers
$s_1,s_2,...,s_L$, where $s_1 = 0$, $s_2 = 1$, and for 
$i \geq 3$,
$$
s_i\ =\ s_j \circ s_k,\ {\rm where\ } j, k < i,
$$
where $\circ$ is either addition, subtraction, or multiplication, and 
where $s_L = n$.  If one expresses such a computation as a string 
indices $(j,k)$ and operations $+, -$ and $\times$, 
then given an integer $N \geq 2$, one can compute $n \pmod{N}$ using only 
$O(L (\log L)(\log N)(\log\log N)^2)$ 
bit operations, by just computing the sequence $s_1,...,s_L$ modulo $N$.
The factor $(\log N)(\log\log N)^2$ in this big-O appears because a product of 
integers modulo $N$ can be computed using Fast Fourier Transforms using
only $O((\log N)(\log\log N)^2)$ bit operations.

We will say that an integer $n$ has {\it computational complexity} $L$
if and only if there exists a computation of lenth $L$ for computing $n$.
\bigskip

In the next section we will prove a general result, which we will use to
prove the following theorem:

\begin{theorem} \label{complexity_theorem}
Given a positive integer $n$ having computational complexity $L$, an 
approximation $A$
to the number of base-$b$ digits of $n$ (see the input specs below), an integer $m \geq 0$,
and a level $y$, there exists an algorithm for computing a level-$y$ approximation
to 
$$
\nu\ =\ \left \{ {n \over b^m} \right \}.
$$ 
This algorithm requires only
$$
O(yL (\log^{2/3} n) \log^C(b + m + y + L + \log n))\ {\rm bits\ of\ memory},
$$
and
$$
O(yL(\log n) \log^C(b+ m + y + L + \log n))\ {\rm bit\ operations.}
$$

The input and output requirements of this algorithm are as follows:

\begin{quote}
Input:  The integers $m$, $y$, and a string representing the
length-$L$ computation of $n$.  Also, the algorithm requires as input an 
integer $A$, which is an approximation to the number $d$ of base-$b$ digits 
of $n$.  This approximation need only satisfy
$$
{1 \over 2}\ <\ {A \over d}\ <\ 2.
$$  
We further restrict $m$ so that $m < 2A+y+1$, since otherwise $\gamma=0$ satisfies
the conclusion of our Theorem.
\end{quote}

\begin{quote}
Output:  The level-$y$ approximation to $\nu$, encoded as a string of 
$y+1$ base-$b$ digits.
\end{quote}
 
\end{theorem}

This theorem requires a little more explanation.  First of all, the input to 
the algorithm will be a string of $O(L \log L + \log m + \log y)$ bits,
which is smaller than the space requirement listed above (when $C > 1$).
The $L\log L$, $\log m$, and $\log y$ terms here account for the number of bits needed to 
specify the length $L$ computation of $n$, the index $m$, and the level $y$, 
respectively.  The approximation $A$ to $d$ requires only $O(\log\log n)$ bits of 
space, and this turns out to be $O(L\log L)$:  To see why this is so, we note 
that any length $L$ computation produces an integer $n < 2^{2^{L-1}}$, which can
be proved by induction.  It follows then that $\log\log n = O(L)$.

The output of the algorithm will be 
a string of $O(y)$ bits, representing the 
base-$b$ approximation $\gamma$ to $\nu$.  The number $\gamma$ will have 
only $y+1$ base-$b$ digits, and will satisfy
$$
|\gamma - \nu|\ <\ {1 \over b^y}.
$$

Perhaps the most surprising aspect of the above theorem is that the indicated
algorithm can require significantly less memory to find the approximation 
to $\nu$ than it does to write down the number $n$, which will have 
$O(\log n)$ base-$b$ digits.  The following corollary of the above theorem
will make this point clear:

\begin{corollary} \label{power_corollary}
Suppose that $a,b \geq 2$ and $m, t, y \geq 1$ are all integers.  There exists an 
algorithm which computes a level-$y$ base-$b$ 
approximation $\gamma$ to $\nu$, where 
$$ 
\nu\ =\  \left \{ {a^t \over b^m} \right \}.
$$
This algorithm requires only
$$
O( yt^{2/3}\log^C (a + b + y + m + t))\ \ {\rm bits\ of\ memory},
$$
where $C > 0$, and performs 
$$
O( yt \log^C (a+b+y+m+t))\ \ {\rm bit\ operations.}
$$
\end{corollary}
Now, the number of bits needed to write down the number $n=a^t$ is
clearly $O(t \log a)$; and yet, if, say, we take $y=1$, this algorithm
requires only $t^{2/3 + o(1)}$ bits of space.  

This corollary follows since $a^t$ has computational complexity
$$ 
L\ =\ O(\log^2( a + t)).
$$ 
To see this, we note that $a^t$ can be generated by
repeated squaring:  If $t = 2^{t_1} + \cdots + 2^{t_s}$, then 
$a^t = a^{2^{t_1}}\cdots a^{2^{t_s}}$.  These numbers $a^{2^h}$
can be computed by starting with $a$; then squaring to get $a^2$;
then squaring again to get $a^4$; then continuing, this produces the list
$a,a^2,a^4,...,a^{2^h}$ after only $h$ multiplications.
\bigskip

The rest of this paper is organized as follows:  In the next section we will
state the Main Theorem (Theorem \ref{main_theorem}) and then use it to
deduce Theorem \ref{complexity_theorem}. 
In section \ref{proof_main_theorem} we give a proof of the 
Main Theorem.  Finally, in section \ref{A_values_section} we give a proof of 
a proposition (Proposition \ref{A_values}), which is an auxillary result 
needed for the proof of the Main Theorem.
\bigskip

\section{Main Theorem and Proof of Theorem \ref{complexity_theorem}}
\bigskip

Theorem \ref{complexity_theorem} is actually a corollary of a more general result
concerning approximations to $\nu$.  In this section we will state  
this result, which will henceforth be called the Main Theorem, 
and then show how to apply it to prove Theorem 
\ref{complexity_theorem}.  The proof of this theorem, as well as a
brief description of the ideas used to prove it, can be found in 
section \ref{proof_main_theorem}; also, at the end of Section 
\ref{alg2_subsection}, we will give a brief statement on how to parallelize 
the algorithm. 

We suppose that $b \geq 2$ is an integer, which is to be the base used; that
$0 < n/a^t < 1$ is some rational number where $n,a \geq 1$, $t \geq 0$ are integers,
and where $n$ has computational complexity $L$; that $y \geq 1$ is some 
level of precision to be used; and finally, that $\mu \geq 0$ is some integer.  
Then, given any pair of integers $S,T$ satisfying
\begin{equation} \label{ST_requirements}
ST\ >\ {3(\log n + (\mu+y+2)\log b) \over \log a},\ a^S\ >\ T^2,
\end{equation}
we have the following

\begin{theorem}[Main Theorem] \label{main_theorem}
Let 
\begin{equation} \label{nu_definition}
\nu\ =\ \left \{b^\mu {n \over a^t} \right \}.
\end{equation}
There exists an algorithm which computes a level-$y$ 
approximation $\gamma$ to this number $\nu$, where the space 
and time requirements of the algorithm are as follows:

\begin{quote} 
Space: $O(yL(S+T) \log^C(y+L+S+T+a+b+\mu+\log t))$ bits of memory, where 
$C > 0$.
\end{quote}

\begin{quote}
Time:  $O(yL(ST+T^3) \log^C(y+L+S+T+a+b+\mu+\log t))$ bit operations.
\end{quote}
\end{theorem}

Here we give more precise information about the input and 
output specifications of the algorithm:

\begin{quote}
Input:  The positive integers $a, t, \mu,y, S$ and $T$, as well as a string of
$O(L \log L)$ characters representing the length-$L$ computation needed to
produce $n$.
\end{quote}

\begin{quote}
Output:  The algorithm will give an approximation 
$\gamma$ to $\nu$.
This approximation will have $y+1$ base-$b$ digits, and 
will satisfy 
$$
|\gamma - \nu|\ <\ {1 \over b^y}.
$$ 
\end{quote}

To prove Theorem \ref{complexity_theorem}, using this result, we let  
$a = b$, $t = 2A$, and $\mu = t - m$.  We note that this gives
$$
{b^\mu n \over a^t}\ =\ {n \over b^m}.
$$
We also let 
\begin{eqnarray}
S\ &=&\ \left \lfloor \left ({3(\log b)(3A + \mu + y + 2) \over \log a}\right )^{2/3} 
\right \rfloor + 1,\nonumber \\
T\ &=&\ \left \lfloor \left ( {3(\log b) (3A + \mu + y + 2) \over \log a} \right )^{1/3} 
\right \rfloor + 1. \nonumber
\end{eqnarray}
We note that this choice of $S$ and $T$ satisfies (\ref{ST_requirements}).
\footnote{To show this, one needs the fact that $3A \log b > \log n$, 
which follows since $n$ has $\leq 2A$ base-$b$ digits.}  

Now, applying the algorithm described in Theorem \ref{main_theorem} with
the parameters indicated above, we get the same output as described in Theorem
\ref{complexity_theorem}.  The running time and space requirements to run
this algorithm are also as stated in Theorem \ref{complexity_theorem} for our
particular choices of $S$ and $T$.
\bigskip

\section{Proof of Theorem \ref{main_theorem}} \label{proof_main_theorem}
\bigskip

Let $S$ and $T$ be as in (\ref{ST_requirements}), and let $r$ and $k$ be 
integers such that
$$
t\ =\ Sk - r,\ 0 \leq r \leq S-1.
$$
Then, we have that
$$
\alpha\ =\ {n \over a^t}\ =\ 
{n a^r \over a^{Sk}}.
$$

The idea of the proof of Theorem \ref{main_theorem} is to approximate 
$b^\mu \alpha$ (and therefore $\nu$) as follows:
\begin{equation} \label{alpha_approximation}
b^\mu \alpha\ =\ \gamma_1 + \cdots + \gamma_T + E,
\end{equation}
where
$$
\gamma_j\ =\ b^\mu {A_j \over a^S - j},
$$
for some rationals $A_1,...,A_T$, and where
\begin{equation} \label{E_inequality}
|E|\ <\ {1 \over b^{y+2}}\ {\rm for\ }ST\ {\rm sufficiently\ large}.
\end{equation}
Then, we will find approximations $\gamma_1',...,\gamma_T'$ to
$\{\gamma_1\},...,\{\gamma_T\}$.
Now, if the precision of these approximations $\gamma_1',...,\gamma_T'$ is 
high enough, and if we we let $\Sigma$ satisfy
$$
\Sigma\ =\ \{\gamma_1' + \gamma_2' + \cdots + \gamma_T'\}, 
$$
then $\Sigma$ will be an approximation to $\nu$; and,
if we then take $\gamma$ to be the closest number to $\Sigma$ having 
$y+1$ base-$b$ digits, then $\gamma$ will be a 
level-$y$ approximation to $\nu$.

We claim that the approximations $\gamma_1',...,\gamma_T'$ to 
$\{\gamma_1\},...,\{\gamma_T\}$ need only have
$$
w\ =\ y\ +\ \left \lfloor {\log T \over \log b} \right \rfloor\ +\ 3,
$$
base-$b$ digits (and be level-$w$ approximations), in order to
guarantee that $\Sigma$ is a level-$y+1$ approximation to $\nu$.  
Note that this would imply that 
$$
| \{\gamma_j\}\ -\ \gamma_j' |\ <\ {1 \over Tb^{y+2}}.
$$
To see only $w$ base-$b$ digits are needed, we note that if these numbers 
$\gamma_j'$ satisfy this last inequality, then by the triangle inequality,
\begin{eqnarray}
|\gamma - \nu|\ &\leq&\ {1 \over b^{y+1}} + |\Sigma - \nu|\ \leq\ 
{1 \over b^{y+1}} + \sum_{j=1}^T |\gamma_j' - \{\gamma_j\}|\ +\ |E|\nonumber \\
&<&\ {1 \over b^{y+1}} + {1 \over b^{y+2}} + {1 \over b^{y+2}}\ \leq\ {1 \over b^y},
\nonumber
\end{eqnarray}
as claimed.

Let us now find a set of values for $A_1,...,A_T$ which make (\ref{alpha_approximation})
hold:  Using the geometric series identity, we have that
\begin{equation} \label{approx_equation}
\sum_{i=1}^T {A_i \over a^S - i}\ =\ \sum_{j=1}^\infty {B_j \over a^{Sj}},
\end{equation}
where
$$
B_j\ =\ A_1 + A_22^{j-1} + A_3 3^{j-1} + \cdots + A_T T^{j-1}.
$$
We seek values for $A_1,...,A_T$ so that
$$
B_j\ =\ \cases{
0, &for $1 \leq j \leq T,\ j \neq k$; \cr
na^r, &for $j=k$.
}
$$
The following Proposition gives the solution we seek

\begin{proposition} \label{A_values}
We have that
\begin{equation} \label{aj_equation}
A_j\ =\ {na^r \left ( {\rm Coef.\ of\ }x^{k-1}\ {\rm in\ } \prod_{h=1 \atop h \neq j}^T 
x-h \right ) \over \prod_{h=1\atop h \neq j}^T j-h};
\end{equation}
and,
\begin{equation} \label{aj_bound}
|A_j|\ \leq\ nTa^S4^T.
\end{equation}
\end{proposition}
Note that this implies
$$
\gamma_j\ =\ b^\mu {(-1)^{T-j} n a^r \left ( {\rm Coef.\ of\ }x^{k-1}\ {\rm in\ }
\prod_{h=1 \atop h \neq j}^T x-h \right ) \over (j-1)! (T-j)! (a^S - j)}.
$$

From this proposition we deduce that
\begin{eqnarray}
|E|\ &=&\ b^\mu \left | \sum_{j=T+1}^\infty {A_1 + A_2 2^{j-1} + \cdots + A_T T^{j-1} \over a^{Sj}}
\right |\nonumber \\
&\leq&\ b^\mu \sum_{j=T+1}^\infty {(nTa^S4^T)T^j \over a^{Sj}} \nonumber \\
&=&\ {na^S4^T T^{T+2}b^\mu \over a^{S(T+1)}} \sum_{j=0}^\infty {T^j \over a^{Sj}} \nonumber \\
&=&\ {na^S4^T T^{T+2}b^\mu \over a^{ST} (a^S - T)}
\ <\ {nb^\mu \over a^{ST/3}}\ <\ {1 \over b^{y+2}}, \nonumber
\end{eqnarray}
for $ST$ large enough; and so, (\ref{E_inequality}) follows.  

We now have all the ingredients necessary to prove Theorem \ref{main_theorem},
which we will give as the following algorithm:
\bigskip

\subsection{Algorithm 1} \label{alg2_subsection}
\bigskip

The input, output, and requirements of this algorithm are as stated in 
Theorem \ref{main_theorem}.  Here are the steps of the algorithm:
\bigskip

{\bf 1.}  Let
$$
w\ =\ y + \left \lfloor {\log T \over \log b} \right \rfloor + 3.
$$
Note that this choice of $w$ satisfies
$$
{T \over b^w}\ \leq\ {1 \over b^{y+2}}.
$$
\bigskip

{\bf 2.}  Set $\Sigma = 0$, and let $r$, $S$, $T$ and $k$ be as described 
at the beginning of this section.
\bigskip

{\bf 3.}  For $j$ from $1$ to $T$ do steps 4 through 8.
\bigskip

{\bf 4.}  Compute
$$
Q\ \leftarrow\ a^S - j.
$$
\bigskip

{\bf 5.}  Set
$$
v\ \leftarrow\ (j-1)! (T-j)! Q.
$$
This number can be computed using $O((T+S) \log^C (T+S))$ bit operations,
and just as much memory (for some $C > 0$).
\bigskip

{\bf 6.}  Apply Algorithm 2 
(given in the next subsection of the paper) to compute 
$$
H\ \leftarrow\ {\rm Coef.\ of\ }x^{k-1}\ {\rm in\ } \prod_{h=1 \atop h \neq j}^T
x-h.
$$
This step requires $O(T^2 \log^C T)$ bit operations and $O(T \log^C T)$ 
bits of memory.
\bigskip

{\bf 7.}  Compute 
$$
u\ \leftarrow\ (-1)^{T-j} n a^r H \pmod{v},\ 0 \leq u \leq v-1.
$$
Since $n$ has computational complexity $L$, this step requires only
$O(L(S+T)\log^C(L+S+T+a))$ bit operations, and just as much memory.
\bigskip

{\bf 8.}  Find a number $\tau$ having $w+1$ base-$b$ digits satisfying 
$$
|\tau - \phi|\ <\ {1 \over b^w},
$$
where
$$
\phi\ =\  \left \{ b^\mu {u \over v} \right \}.
$$
(So, $\tau$ will be a level-$w$ approximation to $\phi$.)

We note that this number $\tau$ can be easily computed by first letting
$$
u_0\ \equiv\ b^\mu u \pmod{v},\ 0 \leq u_0 \leq v-1,
$$
and then noting that
$$
\left \{ {b^\mu u \over v} \right \}\ =\ {u_0 \over v}\ =\ 
(0.r_{-1}r_{-2}...)_b.
$$ 
Then, by finding the first $w+1$ significant digits of $u_0/v$, and 
letting $\tau = (0.r_{-1}...r_{-w-1})_b$ one see that the above inequalities 
are satisfied.
\bigskip

{\bf 9.}  Set
$$
\Sigma\ \leftarrow\ \{ \Sigma + \tau \}.
$$
We only need to do level-$w$ arithmetic in base-$b$ here.

(If $j < T$, then increment $j$ and loop back to step 4.)
\bigskip

{\bf 10.}  (We assume $j = T$.)  Let $\gamma$ be the number having $y+1$
base-$b$ digits which comes nearest to $\Sigma$, and then OUTPUT $\gamma$.
\bigskip

We note that we can perform the operations in steps 3 through 8, with
different values of $j$, in parallel.  For example, given two processors,
we can assign processor 1 to perform steps 3 through 8, with values of 
$j \leq T/2$, and then assign processor 2 to do the same, but with
$T/2 < j \leq T$.  This would result in an two-fold reduction in the running
time, as long as $\mu$ is sufficiently large.  Of course, the memory 
requirements would double, because each of the two processors would require
their own seperate memories.  

More generally, we have that, given $M$ processors, for $\mu$ sufficiently
large, Algorithm 1 can be computed in parallel, resulting in an $M$-fold
reduction in running time, but an $M$-fold increase in memory requirements.
\bigskip

\subsection{Algorithm 2}  
\bigskip

\begin{quote}
Input:  $T,k,j$.
\end{quote}

\begin{quote}
Output:  Coef. of $x^{k-1}$ in $\prod_{h=1 \atop h \neq j}^T x-h$.
\end{quote}

\begin{quote}
Requirements:  The algorithm performs $O(T^2 \log^C T)$ bit operations
(for some $C>0$), but requires only $O(T \log^C T)$ bits of memory.
\end{quote}
\bigskip

{\bf 1.}  Let $P$ be the least integer such that
$$
\Delta\ =\ \prod_{p \leq P \atop p\ {\rm prime}} p\ \geq\ 2^{T+1}T!
$$
Note:  $P=O(T\log T)$, and can be computed using $O(T \log^D T)$
bit operations (for some $D>0$); and so, we can compute and store $P$
within the time and space requirements listed above for the algorithm.
We also note that 
every coefficient of the polynomial in the output specifications
is less than $\Delta/2$ in absolute value. 
\bigskip

{\bf 2.}  Set $\Sigma = 0$.
\bigskip

{\bf 3.}  For each prime $p \leq P$ do steps 4 through 8. 
\bigskip

{\bf 4.}  Compute the polynomial
$$
f(x)\ \equiv\ \prod_{h=1 \atop h \neq j}^T x-h \pmod{p}.
$$ 

Note:  This polynomial can be stored as a length-$T$ coefficient vector,
and the number of bits required to store such a vector is $O(T \log p)
= O(T \log T)$; also, this polynomial can be computed using $O(T \log^D T)$
bit operations by making use of FFT's and a divide-and-conquer 
strategy for polynomial multiplication.  The divide-and-conquer part of
the algorithm can probably best be described as the following recursive 
procedure:  First, we suppose that $L$ is a set of polynomials to be 
producted together modulo $p$, and Product$(L)$ denotes the procedure for
computing this product.  The pseudocode for this procedure is given as
follows:
\bigskip

If $|L|= 1$ (i.e. $L$ has only one polynomial), 
then 

\hskip0.5in RETURN the contents of $L \pmod{p}$;

Else, if $|L| \geq 2$, say $L = \{f_1,...,f_t\}$, then 

\hskip0.5in RETURN 

\hskip0.75in Product$(\{f_1,...,f_{\lfloor t/2 \rfloor} \})
\cdot$Product$(\{f_{\lfloor t/2\rfloor +1},...,f_t\}) \pmod{p}$ 
\bigskip

Now, using FFT's to perform the polynomial multiplication in this second
step (the `Else' step), we see that if the two polynomials being multiplied 
together have degrees $\ell_1$ and $\ell_2$, respectively, then the multiplication
should take no more than $O((\ell_1 + \ell_2)\log^D (\ell_1 + \ell_2+p))$
bit operations.  Now, if we run Product$(L)$ starting with $L$ consisting
of all linear factors $x-h$, $1 \leq h \leq T$, $h \neq j$, then if 
$T-1$ is a power of $2$, the procedure products together 
$(T-1)/2$ pairs of degree $1$ polynomials; 
$(T-1)/4$ pairs of degree $2$ polynomials; and so on, all the way down to 
two polynomials of degree $(T-1)/2$.  
So, the total number of bit operations required to run this producedure is
$$
\ll\ (\log^D T) \sum_{j \leq (\log T)/\log 2 + 1} {T \over 2^j} 2^j\ 
=\ O( T \log^D T).
$$
The memory requirements (in bits) are likewise of the same order.
\bigskip    

{\bf 5.}  Set
$$
H\ \leftarrow\ {\rm Coef.\ of\ }x^{k-1}\ {\rm in\ } f(x) \pmod{p}.
$$
\bigskip

{\bf 6.}  Set
$$
N\ \leftarrow\ (\Delta/p)^{-1} H \pmod{p},\ {\rm where\ }
0 \leq N \leq p-1.
$$
\bigskip 

{\bf 7.}  Set 
$$
\Sigma\ \leftarrow\ \Sigma + {N \Delta \over p}.
$$
\bigskip

{\bf 8.}  Increment the value of $p$, and return to step 4, unless 
$p > P$, in which case we proceed to step 9.
\bigskip

{\bf 9.}  Let $r$ be the least residue in absolute value of $\Sigma 
\pmod{\Delta}$.
\bigskip

{\bf 10.}  Return the value of $r$, and STOP.
\bigskip

It is relatively easy to see that the algorithm requires no more than the 
indicated space and time requirements.  

The idea behind the algorithm is that we use the Chinese Remainder
Theorem to compute the $x^{k-1}$ coefficient of our polynomial, and the computation
in step 7 is just an ``on the fly'' CRT calculation.  This calculation is based on
the following fact:  If $q_1,...,q_h$ are coprime, and if $a_1,...,a_h$ are
any integers, then if we set 
$$
\Delta'\ =\ \prod_{i=1}^h q_i,
$$
and
$$
\Sigma'\ =\ \sum_{i=1}^h b_i {\Delta' \over q_i},\ {\rm where\ }
b_i \equiv a_i (\Delta'/q_i)^{-1} \pmod{q_i},
$$
then
$$
\Sigma'\ \equiv\ a_i \pmod{q_i},\ {\rm for\ every\ }i=1,2,...,h.
$$

One might guess that the coefficient of our polynomial can be computed using
less resources by using a ``Fourier Series'' method; that is,
$$
{\rm Coef.\ }x^{k-1}\ {\rm in\ } f(x)
\ =\ {1 \over T} \sum_{\ell = 0}^{T-1} e^{-2\pi i \ell (k-1)/T} 
f(e^{2\pi i \ell / T}).
$$
It is not obvious (to me) how to do this without using the special form
of the polynomial $f(x)$:  First of all, we would need to maintain
$\gg T$ digits of precision for each term in the sum, since any
particular coefficient of the polynomial $f(x)$ can have size 
$2^{cT \log T}$, for some $c>0$.  Thus, $\gg T^2$ bit operations would be 
needed to compute each term 
$f(e^{2\pi i  \ell/ T})$.  In total, $\gg T^3$ bit operations would be needed to 
evaluate all the terms in the sum.  If one tries to use FFT's to evaluate 
{\it all} the terms in the sum at the same time, this reduces the running 
time to $O(T^2 \log^D T)$ bit operations; however, the memory requirements then 
increase to $\gg T^2$ bits of storage, which is the amount needed to 
store the all numbers $f(e^{2\pi i \ell/T})$, $0 \leq \ell \leq T$ to 
$\gg T$ bits of precision.  Even if we try a discrete version of this method,
where the polynomials are computed, say, modulo $2^k$ for $k \gg T$, and the 
roots of unity are roots of unity modulo $2^k$, we would run into the same
difficulties. 

\bigskip

\section{Proof of Proposition \ref{A_values}} \label{A_values_section}
\bigskip

The $A_i$'s can be computed by solving the equation

\begin{equation} \label{matrix_equation}
\left ( \matrix{ 1 & 1 & 1 & \cdots & 1 \cr
1 & 2 & 3 & \cdots & T \cr
1 & 2^2 & 3^2 & \cdots & T^2 \cr
\vdots & \vdots & \vdots & \cdots & \vdots \cr
1 & 2^{T-1} & 3^{T-1} & \cdots & T^{T-1}
} \right ) \left ( \matrix{A_1 \cr A_2 \cr A_3 \cr \vdots \cr A_T} \right )
\ =\ \left ( \matrix{0 \cr \vdots \cr 0 \cr na^r \cr 0 \cr \vdots \cr 0 } \right ).
\end{equation}
If we call the matrix on the left-hand-side $M$, then
\begin{equation} \label{aj_value}
A_j\ =\ na^r M^{-1}_{j,k},
\end{equation}
where $M^{-1}_{j,k}$ is the entry in the $j$th row, $k$th column of $M^{-1}$.

We will calculate $M_{j,k}^{-1}$ via polynomial interpolation: 
We have that for any set of ordered pairs
$$
(1,b_1),\ (2,b_2),\ ...,\ (T,b_T),
$$
where $b_1,...,b_T \in {\mathbb C}$, there exists a unique degree $T-1$ polynomial
$f(x) \in {\mathbb C}[x]$ such that
$$
f(i)\ =\ b_i,\ {\rm for\ all\ }i=1,2,...,T;
$$
moreover, if we write
$$
f(x)\ =\ c_T x^{T-1} + c_{T-1} x^{T-2} + \cdots + c_2 x + c_1,
$$
then these coefficients $c_i$ can be calculated in two different ways:
The first way is through basic linear algebra, since
\begin{equation}\label{interpolation_matrix}
\left ( \matrix{
1 & 1 & 1 & \cdots & 1 \cr
1 & 2 & 2^2 & \cdots & 2^{T-1} \cr
1 & 3 & 3^2 & \cdots & 3^{T-1} \cr
\vdots & \vdots & \vdots & \vdots & \vdots \cr
1 & T & T^2 & \cdots & T^{T-1}} \right )
\left ( \matrix{c_1 \cr c_2 \cr c_3 \cr \vdots \cr c_T} \right )
\ =\ \left ( \matrix{b_1 \cr b_2 \cr b_3 \cr \vdots \cr b_T} \right ).
\end{equation}
We notice that the matrix on the left-hand-side is $M'$, the transpose of our
matrix $M$.

The second way of calculating the $c_i$'s is by Lagrange interpolation, which gives
\begin{equation} \label{lagrange_equation}
f(x)\ =\ \sum_{i=1}^T b_i \prod_{h=1 \atop h \neq i}^T {x-h \over i-h}.
\end{equation}

Now, if we suppose that
$$
b_i\ =\ \cases{ 0, &if $i \neq j$, \cr
1, &if $i=j$,}
$$
then for this choice of $b_i$'s, one sees from (\ref{interpolation_matrix}) that 
$$
c_k\ =\ (M')^{-1}_{k,j}\ =\ M^{-1}_{j,k}. 
$$
On the other hand, from (\ref{lagrange_equation}) we see that
$$
c_k\ =\ {\rm Coef.\ of\ }x^{k-1}\ {\rm in\ } \prod_{h=1\atop h \neq j}^T {x-h \over j-h}.
$$
Thus,
$$
M^{-1}_{j,k}\ =\ {\rm Coef.\ of\ }x^{k-1}\ {\rm in\ } \prod_{h=1 \atop h \neq j}^T
{x-h \over j-h},
$$
and we conclude from this and (\ref{aj_value}) that (\ref{aj_equation}) holds.

Finally, to prove (\ref{aj_bound}), we note that 
$$
\left | \prod_{h=1 \atop h \neq j}^T j-h \right |\ =\ (T-j)! (j-1)!,
$$
The coefficient of $x^{k-1}$ in the above polynomial is clearly less than
$$
T! {T \choose k-1}\ <\ T! 2^T.
$$
So,
$$
|A_j|\ \leq\ na^S2^T {T! \over (T-j)! (j-1)!}\ =\ jna^S2^T {T \choose j}\ <\ nT a^S 4^T,
$$
which proves (\ref{aj_bound}).  
\bigskip

\section{Acknowledgements}
\bigskip

I would like to thank Richard Hudson for an email he sent to me, which got me interested
in these digit calculation questions, which eventually lead me to prove the theorems
listed above.  I would also like to thank Kevin Hare for pointing out to me that my
algorithm above is highly parallelizable.

\end{document}